\documentclass[twoside,bezier,reqno]{amsart}

\usepackage{epsfig}
\usepackage{amsmath,amsthm,amsfonts,amssymb,amscd,euscript}

\input{epsf}
\newcommand{\filebegin}{\begin{document}}
\newcommand{\fileend}{\end{document}}
\makeatother
\def\thefootnote{}
\newcommand{\lo}{\longrightarrow}
\newcommand{\NMM}{\hspace*{2mm}}

\renewcommand{\baselinestretch}{1.1}

\input{epsf}
\renewcommand{\baselinestretch}{1.1}
\def\n{\noindent}%

\numberwithin{equation}{section}
\def\mapdown#1{\Big\downarrow\rlap
{$\vcenter{\hbox{$\scriptstyle#1$}}$}}
\newtheorem{theorem}{Theorem}[section]
\newtheorem{lemma}[theorem]{Lemma}
\newtheorem{proposition}[theorem]{Proposition}
\newtheorem{corollary}[theorem]{Corollary}

\theoremstyle{definition}
\newtheorem{definition}[theorem]{Definition}
\newtheorem{example}[theorem]{\sc Example}
\newtheorem{xca}[theorem]{Exercise}
\theoremstyle{remark}
\newtheorem{remark}[theorem]{Remark}

\begin{document}

\setcounter{page}{1} \noindent

\vspace*{2cm}
\begin{center}
{\bf\large ENCOUNTER TO THE ADDITIVE PROPERTY OF LOCAL FUNCTION}
 \\[0.5cm]
{{\bf Kulchhum Khatun, Shyamapada Modak and \\Monoj Kumar Das}
 \\[2mm]
Department of  Mathematics, University of Gour Banga\\Mokdumpur, Malda 732103, West Bengal, India \\[2mm]
{\tt E-mail: kulchhumkhatun123@gmail.com }\\
{\tt E-mail: spmodak2000@yahoo.co.in}\\
{\tt E-mail: dmonojkr1@gmail.com}
} \\[2mm]
\end{center}%
\vspace*{0.5cm}
\begin{quotation}
\noindent
{\footnotesize
{\sc Abstract.}
Through this paper we will modify some of the results of \cite{A2021}, \cite{BTP2023}, \cite{KN2010}, \cite{PNB2020}, \cite{TB2017}, \cite{YK2019}, \cite{YK2021} and
 consequently give the modified results.}
\end{quotation}
\\

\n \textbf{2020 Mathematics subject classification: } 54C60, 54G20.

\markboth
{K. Khatun, S. Modak, M. K. Das }
 {Encounter to the Additive Property of Local Function}



\section{Introduction}

Let $X$ be a set and $\wp(X)$ denotes the power set of $X$. An operator $\rho: \wp(X) \rightarrow \wp(X)$ is said be an additive operator if $\rho(A\cup B) = \rho(A) \cup \rho(B)$. Closure operator in a metric space or a topological space is an example of an additive operator. The set of all limit points operator (or derived set operator) $\prime$ is also an example of an additive operator. The local function \cite{JH1990,K1933, V1944}  in an ideal topological space is an another example of an additive operator. The authors, Yalaz and  Kaymakci in \cite{YK2019} and in \cite{YK2021}, Power et al. in \cite{PNB2020}, Beulah et al. in \cite{BTP2023}, Khan and Noiri in \cite{KN2010}, Abbas in \cite{A2021} and Theodore and Beulah in \cite{TB2017} placed the additive property and its consequences wrongly. Through this paper we present the modification of the additive property and its consequences considered by the above authors.

\section{Preliminaries}
Concise discussion about the necessary facts are:

A subcollection $\mathbb{I}$ of $\wp(X)$ is called an ideal on $X$ if $\mathbb{I}$ is closed under heredity and finite additivity property. The local function of a set $A$ in an ideal topological space \cite{JD1999} $(X, \tau, \mathbb{I})$ ( or simply $\mathbb{I}_X$) is, $A^{\star} = \{z\in X| \; U\cap A \notin \mathbb{I},\; U \in \tau(z) \}$, where $\tau(z)$ denotes the set of all open sets containing $z$. Its associated set-valued set function \cite{SELIM2021} is $\psi$ \cite{TN1986} and it is defined as, $\psi(A) = X\setminus (X\setminus A)^{\star}$. A good number of mathematicians have modified the definition of local function by replacing the open set with generalized open sets and open set related closure operators. Then they try to compare the properties of local function with the new types of local functions. One of the property is additive property. But some authors \cite{A2021, BTP2023, KN2010, PNB2020, TB2017, YK2019, YK2021} have  considered these wrongly. The modified local function related to the generalized open sets are as follows:

 A subset $A$ of a topological space $X$ is said to semi-open \cite{H1975, L1963} (resp. preopen \cite{MEE}, $b$-open \cite{DA}, $\beta$-open \cite{EEM} (or semi-preopen \cite{AD}) if $A\subseteq Cl(Int(A))$ (resp. $A\subseteq Int(Cl(A))$, $A\subseteq Int(Cl(A)) \cup Cl(Int(A))$, $A\subseteq Cl (Int(Cl(A)))$), where `Int' and `Cl' stand for the `Interior' and `Closure' operators respectively. $\tau(z)$ is the collection of all open sets containing $z$ whereas $SO(z, X)$ stands for the collection of all semi-open sets containing $z$. Again $PO(z, X)$ (resp. $BO(z, X)$, $\beta O(z, X)$) denotes for the collection of all preopen (resp. $b$-open, $\beta$-open) sets containing $z$. The `semi closure', `pre closure', `$b$ closure' and `$\beta$ closure' operators will be denoted as `$sCl$', `$pCl$' `$bCl$' and `$\beta Cl$' respectively. They are defined in similar fashion with closure operator in the topological space.  Local functions related to these operators are,

 $A^{\star s}$ \cite{KN2010, SS2012} (resp. $A^{\star p}$ \cite{A2021, MM2012}, $A^{\star b}$ \cite{TB2017}, $A^{\star \beta}$ \cite{M2011, TB2017}) $ = \{z \in X|\; U\cap A \notin \mathbb{I} \}$, where $U \in SO(z, X)$ (resp. $U \in PO(z, X)$, $U \in BO(z, X)$, $U \in \beta O(z, X)$). Closure operator related local function is,  $\Gamma (A)$ \cite{AN2013} (resp. $\gamma (A)$ \cite{IM2018}, $\xi_s(A)$ \cite{YK2019}, $\xi_p(A)$, $\xi_b(A)$, $\xi_{\beta}(A)$ \cite{BTP2023}) $ = \{z\in X |\; Cl(U^o) \cap A \notin \mathbb{I} \}$ (resp. $\{z\in X |\; sCl(U^o) \cap A \notin \mathbb{I} \}$, $\{z\in X |\; sCl(U^s) \cap A \notin \mathbb{I} \}$, $\{z\in X| \; pCl(U^p) \cap A \notin \mathbb{I} \}$, $\{z\in X| \; bCl(U^b) \cap A \notin \mathbb{I} \}$, $\{z\in X |\; \beta Cl(U^{\beta}) \cap A \notin \mathbb{I} \}$), where $U^o \in \tau(z)$, $U^s \in SO(z, X)$, $U^p \in PO(z, X)$, $U^b \in BO(z, X)$, $U^{\beta} \in \beta O(z, X)$.

 The associated set-valued set functions \cite{MS2021} of the above operators are as follows: $\psi_{\xi_s}(A) $\cite{YK2021}$ = X\setminus \xi_s(X\setminus A)$,  $\psi_{\xi_p}(A) = X\setminus \xi_p(X\setminus A)$,  $\psi_{\xi_b}(A) = X\setminus \xi_b(X\setminus A)$,  $\psi_{\xi_{\beta}}(A) $\cite{BTP2023}$ = X\setminus \xi_{\beta}(X\setminus A)$, $\psi_{\Gamma}(A) = X\setminus \Gamma(X\setminus A)$ \cite{AN2013}.

\section{Additive property of local function}
Study of limit points of a set in a topological space via generalized open sets may be comprised by the topological ideal. This study has been introduced in \cite{A2021, AN2013, BTP2023, IM2018, KN2010, MM2012, M2011, TB2017}. Whereas limit point studies via ideal have been introduced by Vaidyanathswami \cite{V1944}, Kuratowski \cite{K1933} and some others (see \cite{TRD1990, HH1976, EH1964, JH1990, 2012, SC, MN2019}).

In some research papers recently published in different journals, the authors have considered limit points via generalized open sets and their related closure operators with an ideal (these are also called local functions). These local functions as like similar with the earlier local function \cite{JH1990, K1933, V1944}.  But some properties of these local functions are not same as the kuratowski's local function. One of them is additive property which has been wrongly considered in \cite{A2021, BTP2023, KN2010, PNB2020, TB2017, YK2019, YK2021}.

As a part of the Kuratowski's local function $Cl^{\star}(A) = A\cup A^{\star}$ makes a closure operator and it as like similar with the operator $Cl(A) = A \cup A^{\prime}$. But is not meaning that $A^*= A^{\prime}$. For this, we consider the usual topology on the real line $\mathbb{R}$ and ideal $\mathbb{I}=\{\emptyset, \{0\}\}$. Then for $A = \{\frac{1}{n}|\;  n\in \mathbb{Z}^+\}$, $A^{\prime} = \{0\}$ but $1\in A^{\star}$.

However the operator $()^{\star}$ is not countable additive and it is supported by the following example:

Consider the usual topology $\mathbb{R}_u$ on the real line $\mathbb{R}$. $Q$ denotes the set of  rationals and put $\mathbb{I} = \{\emptyset\}$. Then $\bigcup \limits_{r\in Q}(Cl(\{r\}))= Q$, but $Q^{\star}=Cl(Q)=\mathbb{R}$.

This example has already been cited in \cite{MS2018}.

Now we are recalling the additive property of local function:

\begin{lemma}\cite{JH1990, K1933}
\label{L0}
Let $\mathbb{I}_X$ be an ideal topological space. Then for $A, B\in \wp (X)$, $(A\cup B)^{\star}= A^{\star}\cup B^{\star}$.
\end{lemma}

Crucial steps for proving the Lemma are: Let $x\in (A\cup B)^{\star}$, but $x\notin A^{\star}\cup B^{\star} \implies \exists$ open sets $U_x, V_x\in \tau(x):$ $U_x\cap A\in \mathbb{I}$ and $V_x\cap B\in \mathbb{I} \implies (U_x\cap A)\cup (V_x\cap B)\in \mathbb{I}$ and $(A\cup B)\cap (U_x\cap V_x)\subseteq (U_x\cap A)\cup (V_x\cap B)\in \mathbb{I}\implies x\notin (A\cup B)^{\star}$.

Later, we shall refer these steps. It is noted that the operator $`Cl$' induces a topology on $X$ like the same behaviour of the closure operator $Cl(A)= A\cup A^{\prime}$. This topology is mentioned here by $*$-topology. We have learnt that the property $(A\cup B)^{\star}= A^{\star}\cup B^{\star}$ makes an important role for making the new topology on $X$.

In \cite{KN2010}, the authors Khan and Noiri, in \cite{A2021}, the author Abbas, in \cite{PNB2020}, the authors Powar et al., in \cite{BTP2023}, the authors Beulah et al.,   in \cite{TB2017}, the authors Theodore and Beulah have considered the additive property of the following local functions via different types of generalized open sets:

\begin{lemma}
  Let $\mathbb{I}_X$ be an ideal topological space. Then for $A,B\subseteq X$,
  \begin{enumerate}
    \item $(A\cup B)^{\star s}= A^{\star s}\cup B^{\star s}$
    \item $(A\cup B)^{\star p}= A^{\star p}\cup B^{\star p}$
    \item $(A\cup B)^{\star \beta}= A^{\star \beta}\cup B^{\star \beta}$.
  \end{enumerate}
\end{lemma}

But the above results are not true in general and these have been justified by the following examples:

\begin{example}
\label{sp}
  \begin{enumerate}
    \item Let $X=\{\varpi_1, \varpi_2, \varpi_3, \varpi_4\}$,
     $\tau = \{\emptyset, X, \{\varpi_1\}, \{\varpi_2\},\\ \{\varpi_1, \varpi_2\}\}$, $\mathbb{I}= \{\emptyset, \{\varpi_3\}\}$. Then $SO(X)= \{\emptyset, X, \{\varpi_1\}, \{\varpi_2\}, \{\varpi_1, \varpi_2\},\\ \{\varpi_1, \varpi_3\},  \{\varpi_2, \varpi_3\}, \{\varpi_1, \varpi_4\},\{\varpi_2, \varpi_4\},\{\varpi_1, \varpi_2, \varpi_3\}, \{\varpi_1, \varpi_2, \varpi_4\},\\ \{\varpi_1, \varpi_3, \varpi_4\}, \{\varpi_2, \varpi_3, \varpi_4\}\}$. Then for $E=\{\varpi_1, \varpi_3\}$, $F= \{\varpi_2, \varpi_3\}$, $(E\cup F)^{\star s}\neq E^{\star s}\cup F^{\star s}$.

    \item Let $X=\{\varpi_1, \varpi_2, \varpi_3, \varpi_4\}$, $\tau = \{\emptyset, X, \{\varpi_3, \varpi_4\}, \{\varpi_2, \varpi_3, \varpi_4\}, \\\{\varpi_1, \varpi_3, \varpi_4\}\}$, $\mathbb{I}= \{\emptyset, \{\varpi_1\}\}$. Then $PO(X)= \{\emptyset, X, \{\varpi_3\}, \{\varpi_4\},\\ \{\varpi_3, \varpi_4\}, \{\varpi_1, \varpi_3\}, \{\varpi_2, \varpi_3\}, \{\varpi_2, \varpi_4\}, \{\varpi_1, \varpi_4\},\{\varpi_1, \varpi_2, \varpi_3\},\\ \{\varpi_1, \varpi_3, \varpi_4\}, \{\varpi_1, \varpi_2, \varpi_4\}, \{\varpi_2, \varpi_3, \varpi_4\}\}$ and $\beta O(X)= \{\emptyset, X, \{\varpi_3\},\\ \{\varpi_4\}, \{\varpi_3, \varpi_4\}, \{\varpi_1, \varpi_3\},  \{\varpi_2, \varpi_3\}, \{\varpi_1, \varpi_4\},\{\varpi_2, \varpi_4\},\{\varpi_1, \varpi_2, \varpi_3\},\\ \{\varpi_2, \varpi_3, \varpi_4\}, \{\varpi_1, \varpi_3, \varpi_4\}, \{\varpi_1, \varpi_2, \varpi_4\}\}$. Then for $A=\{\varpi_1, \varpi_3\}$, $B= \{\varpi_1, \varpi_4\}$, $(A\cup B)^{\star p}$ (resp. $(A\cup B)^{\star \beta}$)$\neq A^{\star p}\cup B^{\star p}$ (resp. $A^{\star \beta}\cup B^{\star \beta}$).
  \end{enumerate}
\end{example}

These examples have already been considered in \cite{SS2012}, \cite{MM2012} and \cite{M2011}.
But it is mentioned that, the proof of the Lemma \ref{L0} contains the step, intersection of two open sets is open. However, intersection of two semi-open (resp. preopen, $\beta$-open) sets is not again semi-open (resp. preopen, $\beta$-open) in general.

Consequently, following are wrong results in the papers \cite{PNB2020}, \cite{TB2017}, \cite{BTP2023} and \cite{A2021}:
Since $\beta$-local function and pre-local function are not obey additive property, thus these local functions do not give the closure operator like Kuratowski's local function, $Cl^{\star}(A) = A\cup A^{\star}$. But the authors consider following wrong result :

\begin{lemma}
  Let $\mathbb{I}_X$ be an ideal topological space. Then the operators $Cl^{\star \beta},\; Cl^{\star p}: \wp(X)\rightarrow \wp(X)$ defined by $Cl^{\star \beta}(A)= A\cup A^{\star \beta}$ and $Cl^{\star p}(A)= A\cup A^{\star p}$ are closure operators.
\end{lemma}

Since the operators $( )^{\star p}$ and $( )^{\star \beta}$ are not additive in general, thus the operators in the above lemma are not the closure operator. Therefore they do not induce topologies but in \cite{PNB2020}, \cite{TB2017}, \cite{BTP2023} and \cite{A2021} the authors added the topologies $\tau^{\star p}$ and $\tau^{\star \beta}$ wrongly.

Using the fact of the additive property of pre-local function and $\beta$-local function, the authors in \cite{A2021} and \cite{KN2010} have stated the following result:

\begin{lemma}
Let $\mathbb{I}_X$ be an ideal topological space. Then for $A, B \subseteq X$,
\begin{enumerate}
  \item $A^{\star s} \setminus B^{\star s} = (A\setminus B)^{\star s}\setminus B^{\star s}$.
  \item $A^{\star p} \setminus B^{\star p} = (A\setminus B)^{\star p}\setminus B^{\star p}$.
\end{enumerate}
\end{lemma}

This is also an incorrect result and it will be encountered by the following example:

\begin{example}
\label{E1 }
Let $X = \{\varpi_1, \varpi_2, \varpi_3, \varpi_4 \},\; \tau = \{\emptyset, X, \{\varpi_1 \}, \{\varpi_2 \}, \{\varpi_1, \varpi_2 \} \},\; \mathbb{I} = \{\emptyset, \{\varpi_3 \} \}$. Suppose $A = \{\varpi_1, \varpi_2, \varpi_3 \}$ and $B = \{\varpi_1, \varpi_3 \}$. Then $(A\setminus B)^{\star s} \setminus B^{\star s} = \{\varpi_2 \} \neq \{\varpi_2, \varpi_3, \varpi_4 \} = A^{\star s} \setminus B^{\star s}$.
\end{example}

\begin{example}
Let $X = \{\varpi_1, \varpi_2, \varpi_3, \varpi_4 \},\; \tau = \{\emptyset, X, \{\varpi_3, \varpi_4 \}, \{\varpi_2, \varpi_3, \varpi_4 \}, \\ \{\varpi_1, \varpi_3, \varpi_4 \} \},\; \mathbb{I} = \{\emptyset, \{\varpi_1 \} \}$. Suppose $A = \{\varpi_1, \varpi_3, \varpi_4 \}$ and $B = \{\varpi_1, \varpi_3 \}$. Then $(A\setminus B)^{\star p} \setminus B^{\star p} = \{\varpi_4 \} \neq \{\varpi_1, \varpi_2, \varpi_4 \} = A^{\star p} \setminus B^{\star p}$.
\end{example}

The author Abbas in \cite{A2021} has define associated set-valued set function $\psi_p$\cite{MS2021} by,
$\psi_p(A) = X\setminus (X\setminus A)^{\star p}$ and considered the properties $(i)$ $\psi_p(A\cap B) = \psi_p(A) \cap \psi_p(B)$ and $(ii)$  $\psi_p(A\cup B) = \psi_p(A) \cup \psi_p(B)$. But these properties are false and counter example is,

\begin{example}
\label{E2 }
Let $X = \{\varpi_1, \varpi_2, \varpi_3, \varpi_4 \},\; \tau = \{\emptyset, X, \{\varpi_3, \varpi_4 \}, \{\varpi_2, \varpi_3, \varpi_4 \}, \\ \{\varpi_1, \varpi_3, \varpi_4 \} \},\; \mathbb{I} = \{\emptyset, \{\varpi_1 \} \}$. Suppose $A = \{\varpi_2, \varpi_4 \},\; B = \{\varpi_2, \varpi_3 \},\; E = \{\varpi_3 \}$ and $F = \{\varpi_1, \varpi_2 \}$. Then $\psi_p(A\cap B) \neq \psi_p(A) \cap \psi_p(B)$ and $\psi_p(E\cup F) \neq \psi_p(E) \cup \psi_p(F)$.
\end{example}

The author Abbas in \cite{A2021} used the relation $(A\cap B)^{* p}= A^{* p}\cap B^{* p}$ for the proof of the result $\psi_p(A\cup B) = \psi_p(A) \cup \psi_p(B)$. But the relation $(A\cap B)^{* p}= A^{* p}\cap B^{* p}$ is not true in general.

In Theorem 5.2 of \cite{A2021}, the authors considered the following result:
\begin{theorem}
\label{T1}
Let $\mathbb{I}_X$ be an ideal topological space. Then $\eta=\{A\subseteq X |\; A\subseteq \psi_p(A)\}$ is topology on $X$.
\end{theorem}

But this is also a false result, because the author used the step $A, B \in \eta \Longrightarrow A\cap B \in \eta$. But the step is not correct and it is followed by the following example:

\begin{example}
\label{E3}
Let $X = \{\varpi_1, \varpi_2, \varpi_3, \varpi_4 \},\; \tau = \{\emptyset, X, \{\varpi_3, \varpi_4 \}, \{\varpi_2, \varpi_3, \varpi_4 \}, \\ \{\varpi_1, \varpi_3, \varpi_4 \} \},\; \mathbb{I} = \{\emptyset, \{\varpi_1 \} \}$. Then $\{\varpi_2, \varpi_4 \},\; \{\varpi_2, \varpi_3 \} \in \eta$ but $\{\varpi_2 \} \notin \eta$.
\end{example}

\section{Additive property of closure local function}

The study of closure local function has been introduced Al-Omari and Noiri in \cite{AN2013} and Islam and Modak in \cite{IM2018}. But extension study of the local function has been considered by the authors, Yalaz and Keskin in \cite{YK2019, YK2021} and Beulah et al. in \cite{BTP2023}. In this extension study, some results are not correctly placed in their papers.

Yalaz and Keskin in \cite{YK2019} and Beulah et al. in \cite{BTP2023} discussed the following additive property:

\begin{theorem}
\label{T2}
Let $\mathbb{I}_X$ be an ideal topological space. Then for $A, B \subseteq X$,
\begin{enumerate}
  \item $\xi_s(A\cup B) = \xi_s(A) \cup \xi_s(B)$.
  \item $\xi_{\beta} (A\cup B) = \xi_{\beta} (A) \cup \xi_{\beta}(B)$.
\end{enumerate}
\end{theorem}

But this property is not true. For the proof of this property the authors in \cite{BTP2023} used that `intersection between two $\beta$-open sets is a $\beta$-open set' which is not true.

Counter examples:

\begin{example}
\label{E4}
Let $X = \{\varpi_1, \varpi_2, \varpi_3, \varpi_4 \},\; \tau = \{\emptyset, X, \{\varpi_1 \}, \{\varpi_2 \}, \{\varpi_1, \varpi_2 \} \},\; \\ \mathbb{I} = \{\emptyset, \{\varpi_3 \} \}$. Then semi-open sets containing $\varpi_1$ are: $X,\;\{\varpi_1 \},\; \{\varpi_1, \varpi_2 \},\; \\ \{\varpi_1, \varpi_3 \},\; \{\varpi_1, \varpi_4 \}, \; \{\varpi_1, \varpi_2, \varpi_3 \},\; \{\varpi_1, \varpi_2, \varpi_4 \}, \; \{\varpi_1, \varpi_3, \varpi_4 \}$; semi-closure of the semi-open sets containing $\varpi_1$ are: $\{\varpi_1 \}, X, \; \{\varpi_1, \varpi_3 \},\; \{\varpi_1, \varpi_4 \}, \;\\ \{\varpi_1, \varpi_3, \varpi_4 \}$; semi-open sets containing $\varpi_2$ are: $X,\; \{\varpi_2 \}\; \{\varpi_1, \varpi_2 \}\; \{\varpi_2, \varpi_3 \},\; \\ \{\varpi_2, \varpi_4 \},\; \{\varpi_1, \varpi_2, \varpi_3 \},\; \{\varpi_1, \varpi_2, \varpi_4 \}\; \{\varpi_2, \varpi_3, \varpi_4 \}$; semi-closure of the semi-open sets containing $\varpi_2$ are: $X,\; \{\varpi_2 \},\; \{\varpi_2, \varpi_3 \},\; \{\varpi_2, \varpi_4 \},\; \{\varpi_2, \varpi_3, \varpi_4 \}$; semi-open sets containing $\varpi_3$ are: $X, \; \{\varpi_1, \varpi_3 \},\; \{\varpi_2, \varpi_3 \},\; \{\varpi_1, \varpi_2, \varpi_3 \},\\\; \{\varpi_1, \varpi_3, \varpi_4 \},\; \{\varpi_2, \varpi_3, \varpi_4 \}$; semi-closure of the semi-open sets containing $\varpi_3$ are: $X, \; \{\varpi_1, \varpi_3 \},\;  \{\varpi_2, \varpi_3 \},\; \{\varpi_1, \varpi_3, \varpi_4 \},\; \{\varpi_2, \varpi_3, \varpi_4 \}$; semi-open sets containing $\varpi_4$ are: $X, \; \{\varpi_1, \varpi_4 \},\; \{\varpi_2, \varpi_4 \},\; \{\varpi_1, \varpi_2, \varpi_4 \},\; \{\varpi_1, \varpi_3, \varpi_4 \},\\\; \{\varpi_2, \varpi_3, \varpi_4 \}$; semi-closure of the semi-open sets containing $\varpi_4$ are:  $X,\; \{\varpi_1, \varpi_4 \},\\\; \{\varpi_2, \varpi_4 \},\; \{\varpi_1, \varpi_3, \varpi_4 \},\; \{\varpi_2, \varpi_3, \varpi_4 \}$. Suppose $A = \{\varpi_1, \varpi_3 \},\; B = \{\varpi_2, \varpi_3 \}$, then $\xi_s(A) = \{\varpi_1 \},\; \xi_s(B) = \{\varpi_2\},\; \xi_s(A\cup B) =  \{\varpi_1, \varpi_2, \varpi_3, \varpi_4 \}$.
\end{example}

\begin{example}
\label{E5}
Let $X = \{\varpi_1, \varpi_2, \varpi_3, \varpi_4 \},\; \tau = \{\emptyset,\; X,\; \{\varpi_3, \varpi_4 \}, \; \{\varpi_2, \varpi_3, \varpi_4 \}, \\ \; \{\varpi_1, \varpi_3, \varpi_4 \} \}, \;  \mathbb{I} = \{\emptyset, \{\varpi_1 \} \}$. Suppose $A = \{\varpi_1, \varpi_3 \},\; B = \{\varpi_1, \varpi_4 \}$. Then $\xi_{\beta}(A) = \{\varpi_3 \},\; \xi_{\beta}(B) = \{\varpi_4 \},\; \xi_{\beta}(A\cup B) =  \{\varpi_1, \varpi_2, \varpi_3, \varpi_4 \}$.
\end{example}

\begin{example}
\label{E6}
Let $X =\{\varpi_1, \varpi_2, \varpi_3, \varpi_4 \},\; \tau = \{\emptyset,\; X,\; \{\varpi_3, \varpi_4 \}, \; \{\varpi_2, \varpi_3, \varpi_4 \}, \\ \; \{\varpi_1,  \varpi_3, \varpi_4 \} \}, \; \mathbb{I} =\{\emptyset, \{\varpi_1 \} \}$. Suppose $E = \{\varpi_1, \varpi_3 \},\; F = \{\varpi_1, \varpi_4 \}$. Then $\xi_p(E) = \{\varpi_3 \},\; \xi_p(F) = \{\varpi_4 \},\; \xi_p(E\cup F) =  \{\varpi_1, \varpi_2, \varpi_3, \varpi_4 \}$.
\end{example}

In \cite{YK2021} and in \cite{BTP2023}, the authors state and prove the following theorems:

\begin{theorem}
\label{000001}
Let $\mathbb{I}_X$ be an ideal topological space. Then for $A, B \subseteq X$,
\begin{enumerate}
  \item $\psi_{\xi_s}(A\cap B) = \psi_{\xi_s}(A) \cap \psi_{\xi_s}(B)$.
  \item $\psi_{\xi_{\beta}}(A\cap B) = \psi_{\xi_{\beta}}(A) \cap \psi_{\xi_{\beta}}(B)$.
\end{enumerate}
\end{theorem}

In the proof of the theorem the authors used the additive property of $\xi_s$ and $\xi_{\beta}$ operators which are not correct. Thus this theorem is not true.

\begin{theorem}
\label{T3}
Let $\mathbb{I}_X$ be an ideal topological space. Then
\begin{enumerate}
  \item $\sigma_{\xi_s} = \{A\subseteq X| \; A\subseteq \psi_{\xi_s}(A) \}$ forms a topology on $X$ \cite{YK2021}.
  \item $\sigma_{\xi_{\beta}} = \{A\subseteq X| \; A\subseteq \psi_{\xi_{\beta}}(A) \}$ forms a topology on $X$ \cite{BTP2023}.
\end{enumerate}
\end{theorem}
 This is also a wrong result, as they used the Theorem \ref{000001} which is wrong. In this purpose the counter examples are:

\begin{example}
\label{E7 }
Let $X = \{\varpi_1, \varpi_2, \varpi_3, \varpi_4 \},\; \tau = \{\emptyset,\; X,\; \{\varpi_1 \},\; \{\varpi_2 \}, \; \{\varpi_1, \varpi_2 \} \}, \; \\ \mathbb{I} = \{\emptyset, \{\varpi_3 \} \}$.
semi-closure of the semi-open sets containing $\varpi_1$ are: $X,\; \{\varpi_1 \},\; \\ \{\varpi_1, \varpi_3 \},\; \{\varpi_1, \varpi_4 \},\; \{\varpi_1, \varpi_3, \varpi_4 \}$; semi-closure of the semi-open sets containing $\varpi_2$ are: $X,\; \{\varpi_2 \},\;\{\varpi_2, \varpi_3 \},\;\{\varpi_2, \varpi_4 \},\;\{\varpi_2, \varpi_3, \varpi_4 \}$; semi-closure of the semi-open sets containing $\varpi_3$ are: $X,\; \{\varpi_1, \varpi_3 \}, \; \{\varpi_2, \varpi_3 \},\;\{\varpi_1, \varpi_3, \varpi_4 \}, \\ \; \{\varpi_2, \varpi_3, \varpi_4 \}$; semi-closure of the semi-open sets containing $\varpi_4$ are: $X,\; \{\varpi_1, \varpi_4 \}, \;\\ \{\varpi_2, \varpi_4 \}, \; \{\varpi_1, \varpi_3, \varpi_4 \}, \; \{\varpi_2, \varpi_3, \varpi_4 \}$. Consider $A = \{\varpi_2, \varpi_4 \},\; B = \{\varpi_1, \varpi_4 \}$. Then $\psi_{\xi_{s}}(A) = X\setminus \xi_{s} (\{\varpi_1, \varpi_3\}) = X \setminus \{\varpi_1 \} = \{\varpi_2, \varpi_3, \varpi_4 \}$ and  $\psi_{\xi_{s}}(B) = X\setminus \xi_{s}(\{\varpi_2, \varpi_3 \}) = X\setminus \{\varpi_2 \} = \{\varpi_1, \varpi_3, \varpi_4 \}$. Now $\psi_{\xi_{s}}(A\cap B) = X\setminus  \xi_{s}(\{\varpi_1, \varpi_2, \varpi_3 \}) = X\setminus X$. Thus we have, $A\subseteq \psi_{\xi_s} (A),\; B\subseteq \psi_{\xi_s} (B)$. But $A\cap B  \nsubseteq \psi_{\xi_{s}} (A\cap B)$.
\end{example}

\begin{example}
\label{E8}
Let $X = \{\varpi_1, \varpi_2, \varpi_3, \varpi_4 \},\; \tau = \{\emptyset,\; X,\; \{\varpi_3, \varpi_4 \}, \; \{\varpi_2, \varpi_3, \varpi_4 \}, \\ \; \{\varpi_1, \varpi_3, \varpi_4 \} \}, \;  \mathbb{I} = \{\emptyset, \{\varpi_1 \} \}$.

$\beta$-closure of the $\beta$-open sets containing $\varpi_1$ are: $X,\; \{\varpi_1, \varpi_3\},\; \{\varpi_1, \varpi_4 \},\;\\ \{\varpi_1, \varpi_2, \varpi_3 \},\; \{\varpi_1, \varpi_2, \varpi_4 \}$; $\beta$-closure of the $\beta$-open sets containing $\varpi_2$ are: $X,\;\{\varpi_2, \varpi_3\},\;\{\varpi_2, \varpi_4 \},\;\{\varpi_1, \varpi_2, \varpi_3\},\; \{\varpi_1, \varpi_2, \varpi_4\}$; $\beta$-closure of the $\beta$-open sets containing $\varpi_3$ are: $X,\; \{\varpi_3 \},\; \{\varpi_1, \varpi_3 \},\;\{\varpi_2, \varpi_3\},\; \{\varpi_1, \varpi_2, \varpi_3 \}$; $\beta$-closure of the $\beta$-open sets containing $\varpi_4$ are: $X,\; \{\varpi_4 \},\;\{\varpi_1, \varpi_4 \},\;\{\varpi_2, \varpi_4\},\;\\ \{\varpi_1, \varpi_2, \varpi_4 \}$.

Let $A = \{\varpi_2, \varpi_4 \},\; \; B = \{\varpi_2, \varpi_3 \}$. Then $\psi_{\xi_{\beta}} (A) = X\setminus \xi_{\beta} (\{\varpi_1, \varpi_3 \}) = X \setminus \{\varpi_3 \} = \{\varpi_1, \varpi_2, \varpi_4 \}$ and $\psi_{\xi_{\beta}} (B) = X\setminus \xi_{\beta} (\{\varpi_1, \varpi_4 \}) = X\setminus \{\varpi_4 \} = \{\varpi_1, \varpi_2, \varpi_3 \}$. Now, $\psi_{\xi_{\beta}} (A\cap B) = X\setminus \xi_{\beta}(\{\varpi_1, \varpi_3, \varpi_4 \}) = X\setminus X$. Thus we have, $A\subseteq \psi_{\xi_{\beta}} (A),\; B\subseteq \psi_{\xi_{\beta}} (B)$. But $A\cap B  \nsubseteq  \psi_{\xi_{\beta}} (A\cap B)$.
\end{example}



\providecommand{\bysame}{\leavevmode\hbox
to3em{\hrulefill}\thinspace}


\end{document}